\documentclass[notitlepage,11pt]{article}
\usepackage{amssymb,amsmath,geometry,
}

\catcode`\@=11
\@addtoreset{equation}{section}

\catcode`\@=12

\allowdisplaybreaks

\usepackage{color}

\def\R{\mathbb{R}}

\def\f{\varphi}

\def\irn{\int\limits_{\R^n}}

\def\eps{\varepsilon}
\def\div{{\rm div}}

\def\Ds{\left(-\Delta\right)^{\!s}\!}  

\def\DsNROmega{(-\Delta_{\Omega}^{\!N}\!\!\left.\right)^s_{{\rm R}}\!}  

\def\DsNSOmega{(-\Delta_{\Omega}^{\!N}\!\!\left.\right)^{s}_{{\rm Sr}}\!}  

\def\DsNSpOmega{(-\Delta_{\Omega}^{\!N}\!\!\left.\right)^{s}_{{\rm Sp}}\!}  

\def\DsSpOmega{(-\Delta_{\Omega}\!\!\left.\right)^{s}_{{\rm Sp}}\!}  

\def\sstar{{2^*_s}}

\def\Hs{\widetilde{H}^s(\Omega)}							
\def\complement{{\Omega^\mathsf{c}}}

\def\Gcomplement{G^\mathsf{c}}

\def\tr{\tilde r}
\def\tk{\tilde k}
\def\tB{\widetilde B}
\def\tw{\widetilde w}

\def\bp{{\bar p}}

\def\proof{\noindent{\textbf{Proof. }}}
\def\QED{\hfill {$\square$}\goodbreak \medskip}

\newtheorem{Theorem}{Theorem}[section]
\newtheorem{Lemma}[Theorem]{Lemma}
\newtheorem{Proposition}[Theorem]{Proposition}
\newtheorem{Corollary}[Theorem]{Corollary}
\newtheorem{Remark}[Theorem]{Remark}

\begin{document}

\title 
{Strong maximum principles for fractional Laplacians}

\author{Roberta Musina\footnote{Dipartimento di Scienze Matematiche, Informatiche e Fisiche, Universit\`a di Udine,
via delle Scienze, 206 -- 33100 Udine, Italy. Email: {roberta.musina@uniud.it}. 
{Partially supported by
Miur-PRIN 2009WRJ3W7-001.}}~ and
Alexander I. Nazarov\footnote{
St.Petersburg Department of Steklov Institute, Fontanka 27, St.Petersburg, 191023, Russia, 
and St.Petersburg State University, 
Universitetskii pr. 28, St.Petersburg, 198504, Russia. E-mail: al.il.nazarov@gmail.com.
Partially supported by RFBR grant 17-01-00678.
}
}

\date{}

\maketitle

\begin{abstract}
We prove strong maximum principles for a large class of nonlocal operators of the order $s\in(0,1)$,
that includes the Dirichlet, the Neumann Restricted (or Regional) and the Neumann Semirestricted
Laplacians.
\footnotesize

\medskip

\medskip

\noindent
\textbf{2010 Mathematics Subject Classfication:} 47A63; 35A23.
\end{abstract}

\normalsize

\bigskip

\maketitle

\section{Introduction}

In this paper we prove strong maximum principles for a large class of fractional Laplacians
of order $s\in(0,1)$, including

\noindent the {\em  Dirichlet Laplacian} $~~~\!\displaystyle{\Ds u(x)=C_{n,s}~\!\cdot~\!{\rm P.\!V.}\irn \frac{u(x)-u(y)}{|x-y|^{n+2s}}~dy}$,\\
the {\em Restricted Neumann Laplacian}
$~~~\!\displaystyle{\DsNROmega u(x)=~\!C_{n,s}~\!\cdot~\!{\rm P.\!V.}\int
\limits_{\Omega} \frac{u(x)-u(y)}{|x-y|^{n+2s}}~dy}$

\noindent {(also called {\em Regional} Laplacian), }and intermediate operators, such as the {\em Semirestricted Neumann Laplacian}
\begin{equation}
\label{eq:SR}
\DsNSOmega u=\chi_{\Omega}\cdot(-\Delta)^{\!s} u+\chi_{\complement}\cdot\DsNROmega u~\!.
\end{equation}
Here {$\Omega$ is a domain} in $\R^n$, $n\ge 1$, $\complement=\R^n\setminus\Omega$, 
$\chi_V$ is the characteristic function of the set $V\subset \R^n$, 
$C_{n,s}=\frac{s2^{2s}\Gamma(\frac{n}{2}+s)}{\pi^{\frac{n}{2}}\Gamma(1-s)}$ and ''P.\!V.'' means ''principal value''.

The standard basic reference for the popular and largely studied operator $\Ds$ is the monograph \cite{Tr}. 
Restricted {Neumann} Laplacians appear as generators of so-called censored processes; a
vaste literature about the  operator $\DsNROmega$ 
is available as well, see for instance \cite{G1, G2} and \cite{W1, W2, W3}, where 
Neumann, Robin and mixed boundary value problems on not necessarily regular domains $\Omega$ are studied.
The Semirestricted {Neumann} Laplacian $\DsNSOmega$ has been proposed 
in \cite{DRoV} to set up an {alternative approach to  Neumann problems}, and can be used to study non-homogeneous Dirichlet problems for $\Ds$,
see for instance the survey paper \cite{Ro}.
\medskip

\medskip

In this paper we propose a unifying approach to handle, in particular, all fractional Laplacians above. Let 
us describe the  class of nonlocal operators 
we are interested in. \medskip

Consider a domain $\Omega\subseteq\R^n$ and open sets $U_1, U_2\subseteq\R^n$ such that $\Omega\subseteq U_1\cap U_2$.
We put 
\begin{equation*}
Z=\big(U_1\times U_2\big)~\!\cup ~\!\big(U_2\times U_1\big)\subseteq \R^n\times\R^n~,
\end{equation*}
so that $\Omega\times\Omega\subseteq Z$.
For $s\in(0,1)$ we introduce the space
$$
X^s(\Omega;Z)=\big\{~\! u:{U_1\cup U_2} 
\to \R~\text{measurable}~\big|~
\frac{u(x)}{1+|x|^{n+2s}}\in L^1(U_1\cup U_2), ~u\in H^s_{\rm loc}(\Omega)~\big\}.
$$
Notice that, in particular, $X^s(\Omega;Z)$ contains functions $u\in L^1_{\rm loc}(U_1\cup U_2)$ such that
\begin{equation}
\label{eq:qf}
\mathcal E_s(u;Z):=\frac{C_{n,s}}{2}\iint\limits_{\!\! Z} \frac{(u(x)-u(y))^2}{|x-y|^{n+2s}}~dxdy
\end{equation}
is finite. 
For $u\in X^s(\Omega;Z)$ we introduce the distribution 
${\mathfrak L}^s_Z u\in \mathcal D'({\Omega})$ defined via
\begin{equation*}
\langle {\mathfrak L}^s_Z u, \f\rangle= \frac{C_{n,s}}{2}\iint\limits_{\!\! Z} \frac{(u(x)-u(y))(\f(x)-\f(y))}{|x-y|^{n+2s}}~dxdy~,
\quad \f\in {\cal C}^\infty_0({\Omega})~\!,
\end{equation*}
{see  Lemma \ref{L:dual}.}

{Notice that the operator ${\mathfrak L}^s_Z $ might be used in modeling symmetric,
(possibly) censored L\'evy flights of a particle that can only jump from points $x\in U_1$ to points $y\in U_2$, and vice versa.}
Our  approach can be easily generalized for a wider class of kernels $\frac{A(x,y)}{|x-y|^{n+2s}}$ with $A$ measurable, symmetric, bounded and bounded away from zero.

We understand the inequality ${\mathfrak L}^s_Zu\ge 0$ in $\Omega$ in distributional sense, that is,
$$
\langle {\mathfrak L}^s_Z u, \f\rangle\ge 0,\qquad\text{if}\quad \f\in  {\cal C}^\infty_0(\Omega),~ \f\ge0. 
$$

In our main result, see Theorem \ref{T:mp}, we provide a strong maximum principle for solutions to ${\mathfrak L}^s_Z u\ge0$ in $\Omega$, with no assumptions on $\Omega$. 
To prove Theorem \ref{T:mp} we follow the outlines of the arguments in \cite[Theorems 2.4 and 2.5]{IMS},
that cover the case ${\mathfrak L}^s_Z=\Ds$, $\Omega$  bounded and smooth, $n\ge 2$, $u\in  H^s(\R^n)$
and $u\ge 0$ in $\R^n\setminus\overline\Omega$. We cite also \cite[Theorem 1.2]{DQ} for a related result involving
the fractional Dirichlet $p$-Laplacian.


\medskip

The paper is organized as follows. In Section~2 we prove some auxiliary statements. Section~3 is devoted to Caccioppoli type estimates
and to De Giorgi{-type maximum estimates for (sub)solutions. 
In Section~4 we state and prove Theorem \ref{T:mp}, and formulate corresponding results for the Dirichlet, Restricted
and Semirestricted Neumann Laplacians.

{In the Appendix we collect some more strong maximum principles 
for nonlocal Laplacians. First, {we formulate a strong maximum principle}
for $\Ds$ that is essentially contained in the remarkable paper \cite{Sil} by Silvestre,
who extended the classical theory
of superhamonic functions to the case of fractional Laplacian.
Then we discuss strong maximum principles for spectral fractional Laplacians. 
{The {\em Spectral Dirichlet Laplacian} $\DsSpOmega$ (also called  {\em the Navier Laplacian}) is widely studied. }
Notice that for $\Omega=\R^n$ we have $\DsSpOmega=\Ds$, for other $\Omega$ these operators differ, see \cite{FL} for some integral and pointwise inequalities 
between them. The {\em Spectral Neumann Laplacian} $\DsNSpOmega$ is less investigated; we limit ourselves to cite \cite{AEW, CaSt, 
G} and references therein.}

\bigskip
{\small 
\noindent
{\bf Notation.} 
Here we recall some basic notions taken from  \cite{Tr}. For $Z\subseteq\R^n\times\R^n$ and $u$ measurable, let $\mathcal E_s(u;Z)$ be the quadratic form in
(\ref{eq:qf}). We put
$$
H^s(\R^n)=\big\{u\in L^2(\R^n)~|~\mathcal E_s(u;\R^n\!\times\!\R^n)<\infty\big\},
$$
that is an Hilbert space with respect to the norm 
$$
\|u\|^2_{H^s(\R^n)}:
=\mathcal E_s(u;\R^n\!\times\!\R^n)+\|u\|^2_{L^2(\R^n)}.
$$
For any domain $G\subset\R^n$, we introduce the following {closed subspace of $H^s(\R^n)$:}
$$
\widetilde H^s(G)=\big\{u\in H^s(\R^n)~|~u=0~~\text{on $\R^n\setminus \overline G$}~\big\}~\!,
$$
and its dual space $\widetilde H^s(G)'$. 

We write $u\in H^s_{\rm loc}(\Omega)$ if for any $G\Subset\Omega$, the function $u$
is the restriction to $G$ of some $v \in H^s(\R^n)$, and we put
$$
\|u\|_{H^s(G)}:=\inf\big\{\|v\|_{H^s(\R^n)}~~|~~v=u~~\text{on $G$}~\big\}.
$$
It is well known that $u\in H^s_{\rm loc}(\Omega)$ if and only if $\eta u\in \widetilde H^s(\Omega)$ 
for any $\eta\in {\cal C}^\infty_0(\Omega)$, see for instance \cite[Subsection 4.4.2]{Tr}.



We adopt the following standard notation:

$B_{r}(x)$ is the Euclidean ball of radius $r$ centered at $x$, and $B_r=B_r(0)$;

$u^\pm=\max\{\pm u, 0\}$; $\sup\limits_G u$ and $\inf\limits_G u$ stand for essential supremum/infimum of the measurable function $u$ on the measurable set $G$;


Through the paper, all constants depending only on $n$ and $s$ are denoted by $c$. To indicate that a constant depends on other quantities 
we list them in parentheses: $c(\dots)$.}

\noindent

\section{Preliminaries}\label{2}

For any {function $\f$ on $\R^n$ 
we put}
\begin{equation}
\label{eq:3,1}
\Psi_{\!\f}(x,y)=\frac{(\f(x)-\f(y))^2}{|x-y|^{n+2s}}.
\end{equation}

\begin{Lemma}
\label{L:dual}
Let $u\in X^s(\Omega;Z)$. Then 
${\mathfrak L}^s_Z u$ is a well defined distribution in $\Omega$. Moreover, 
for any Lipschitz domain $G\Subset \Omega$, we have  
${\mathfrak L}^s_Z u\in \widetilde H^s(G)'$ and
$$\iint\limits_{G\times G}|u(x)u(y)|\Psi_{\!\f}(x,y)\,dxdy<\infty\qquad\text{for any $\f\in {\cal C}^\infty_0(G)$}.
$$
\end{Lemma}

\proof
Let $\f\in {\cal C}^\infty_0(\Omega)$. In order to have that ${\mathfrak L}^s_Z u$ is well defined we need to show that
$$
g(x,y):= \frac{(u(x)-u(y))(\f(x)-\f(y))}{|x-y|^{n+2s}}\in L^1(Z).
$$
Take two Lipschitz domains $G, \widetilde G$, such that
$\text{supp}(\f)\subset G\Subset \widetilde G\Subset \Omega$. From  $u\in H^s(\widetilde G)$, we have 
$$
\iint\limits_{\widetilde G\times\widetilde G}
|g(x,y)|~\!dxdy\le
\|u\|_{H^s(\widetilde G)}\|\f\|_{H^s(\R^n)}\le
c( G)\|u\|_{H^s(\widetilde G)}\|\f\|_{\widetilde H^s(G)}.
$$
Next, since $\f$ vanishes outside $ G$, and since
\begin{equation}
\label{eq:sets}
\big[Z\setminus(\widetilde G\times \widetilde G)\big]\setminus(\Gcomplement\times\Gcomplement)=
\left[ G\times\big((U_1\cup U_2)\setminus \widetilde G\big)\right]~\!\cup~\! 
\left[\big((U_1\cup U_2)\setminus \widetilde G\big)\times  G\right],
\end{equation}
it is enough to prove that $g\in L^1( G\times\big((U_1\cup U_2)\setminus \widetilde G))$. We have
\begin{eqnarray*}
&&\iint\limits_{ G\times((U_1\cup U_2)\setminus \widetilde G)} |g(x,y)|~\!dxdy \le
\int\limits_G|\f(x)|\bigg[\int\limits_{(U_1\cup U_2)\setminus \widetilde G}\frac{|u(x)|+|u(y)|}{|x-y|^{n+2s}}~\!dy\bigg]dx\\
&&\le c({\rm dist}(G,\partial\widetilde G))\Big(\|\f\|_{L^2(G)}\| u\|_{L^2(G)}+\|\f\|_{L^1(G)}\int\limits_{U_1\cup U_2} \frac{|u(y)|}{1+|y|^{n+2s}}~\!dy\Big)
\le c(G,\widetilde G,u)\|\f\|_{\widetilde H^s(G)}.
\end{eqnarray*}
We proved that ${\mathfrak L}^s_Z u\in \mathcal D'(\Omega)$ and actually ${\mathfrak L}^s_Z u\in \widetilde H^s(G)'$,
 by the density of ${\cal C}^\infty_0(G)$ in 
$\widetilde H^s(G)$.} 

Further, take again $\f\in {\cal C}^\infty_0(G)$ and  notice that $\Psi_{\!\f}(x,~\!\cdot~\!)\in L^1(\R^n)$ for any $x\in \R^n$, because
$$
\Psi_{\!\f}(x,y)\le c(\f)\Big(\frac{\chi_{\{|x-y|<1\}}}{|x-y|^{n-2(1-s)}}+\frac{\chi_{\{|x-y|>1\}}}{|x-y|^{n+2s}}\Big).
$$
Actually $\displaystyle{\irn\Psi_{\!\f}(x,y)~\!dy\le c(\f)}$, and 
{by the Cauchy-Bunyakovsky-Schwarz inequality we infer}
$$
\iint\limits_{G\times G}|u(x)u(y)|\Psi_\f~\!dxdy\le \iint\limits_{G\times G}|u(x)|^2\Psi_{\!\f}~\!dxdy\le c(\f)\int\limits_{G}|u(x)|^2~\!dx\le
c(u,\f,G)<\infty~\!.
$$
The lemma is proved.
\QED

Next, for any domain $G\subseteq U_1\cap U_2 $ we introduce the
{\bf relative killing measure} $M^Z_{G}\in L^\infty_{\rm loc}(G)$,
$$
M^Z_{G}(x)={C_{n,s}}\int\limits_{(U_1\cup U_2)\setminus G}\frac{dy}{|x-y|^{n+2s}}~\!,\quad x\in G~\!.
$$
When $U_1\cup U_2=\R^n$, that happens  for instance in the Dirichlet and in the Semirestricted cases,
see Section \ref{S:main},
 the weight
$M^Z_{G}$ coincides with so-called {\em killing measure} 
of the set $G$:
$$
M_{G}(x):=M^{\R^n\!\times\R^n}_{G}\!(x)= {C_{n,s}}\int\limits_{\R^n\setminus G}\frac{dy}{|x-y|^{n+2s}}~\!.
$$
In the Restricted case we have
$U_1\cup U_2=\Omega$ and $M^Z_{G}=M^{\Omega\times\Omega}_G$;  if $G\subset\Omega$ then $M^{\Omega\times\Omega}_G$ is the 
difference between the killing measures of the sets $G$ and $\Omega$. 


\begin{Lemma}
\label{L:1}
Let $G\subseteq U_1\cap U_2$ be a Lipschitz domain. If $u\in\widetilde H^s(G)$, then 
$$\displaystyle{\mathcal E_s(u;Z)=\mathcal E_s(u;G\times G)+ \int\limits_{G}~\! M^Z_{G}(x)|u(x)|^2~\!dx}$$
and in particular $u$ is square integrable on $G$ with respect to the measure $M^Z_G(x)dx$.
\end{Lemma}

\proof
Trivially $\mathcal E_s(u;Z)<\infty$, as $u\in \widetilde H^s(G)\hookrightarrow H^s(\R^n)$.
Since $u$ vanishes on $\overline G^{\,\mathsf{c}}$, using (\ref{eq:sets}) with $\widetilde G=G$ we have 
\begin{eqnarray*}
\mathcal E_s(u;Z)&=&
\mathcal E_s(u;G\times G)+2\mathcal E_s(u;G\times\big[(U_1\cup U_2)\setminus G\big])
\\
&=&\mathcal E_s(u;G\times G)+C_{n,s}\int\limits_{G}|u(x)|^2\,\Big(\int\limits_{(U_1\cup U_2)\setminus G} \frac{dy}{|x-y|^{n+2s}}~\!
\Big)dx,
\end{eqnarray*}
and the lemma is proved. 
\QED

%

The next two elementary lemmata deal with certain quantities, depending on functions $u\in X^s(\Omega;Z)$, that will be involved in the crucial {Caccioppoli-type inequality} in the next section.

For $u\in X^s(\Omega;Z)$ and for any domain $G\subseteq \Omega$ we use Lemma \ref{L:dual} to introduce the distribution
$$
\langle((-\Delta_{(U_1\cup U_2)\setminus G}^{\!N}\!\!\left.\right)^s_{{\rm R}}\! u),\f\rangle
:=\frac{C_{n,s}}{2}\iint\limits_{G\times[(U_1\cup U_2)\setminus G]}\frac{(u(x)-u(y))(\f(x)-\f(y))}{|x-y|^{n+2s}}~\!dxdy,~~\f\in {\cal C}^\infty_0(G),
$$
that is the restriction on $G$ of the Regional Laplacian of $u$ relative to the set $(U_1\cup U_2)\setminus G$.

\begin{Lemma}
\label{L:3}
Let $G\subseteq \Omega$ be a domain, $u\in X^s(\Omega;Z)$. Then for any $\f\in {\cal C}^\infty_0(G)$
\begin{equation}
\label{eq:L3}
\int\limits_G |u(x)||\f(x)|^2\Big(\!\!\int\limits_{(U_1\cup U_2)\setminus G}\!\!\frac{|u(x)-u(y)|}{|x-y|^{n+2s}}~\!dy\Big)dx<\infty.
\end{equation}
In particular, $u\cdot (-\Delta_{(U_1\cup U_2)\setminus G}^{\!N}\!\!\left.\right)^s_{{\rm R}}\! u\in L^1_{\rm loc}(G)$.
\end{Lemma}
 
 \proof 
Similarly as in the proof of Lemma \ref{L:dual}, we estimate the integral in (\ref{eq:L3}) by 
 \begin{multline*}
\int\limits_G|u(x)| |\f(x)|^2\bigg[\int\limits_{(U_1\cup U_2)\setminus  G}\frac{|u(x)|+|u(y)|}{|x-y|^{n+2s}}~\!dy\bigg]dx\\
\le c({\rm dist}(\text{supp}(\f),\partial G))\|\f\|^2_{L^\infty(G)}\Big(\| u\|^2_{L^2(\text{supp}(\f))}+\|u\|_{L^1(\text{supp}(\f))}\int\limits_{U_1\cup U_2} \frac{|u(y)|}{1+|y|^{n+2s}}~\!dy\Big)<\infty~\!,
\end{multline*}
{and the lemma follows.}
\QED

\begin{Lemma}
\label{L:truncation}
If  $u\in X^s(\Omega;Z)$, then $u^\pm\in X^s(\Omega;Z)$; moreover
for any $G\Subset \Omega$ we have $$\mathcal E(u^\pm;G\times G)< \mathcal E(u; G\times G),$$ unless $u$ has constant sign on $G$.
\end{Lemma} 

\proof
We compute
$$
(u(x)-u(y))^2-(u^+(x)-u^+(y))^2=(u^-(x)-u^-(y))^2+2\Big(u^+(x)u^-(y)+u^-(x)u^+(y)\Big)\ge 0.
$$
Thus $\mathcal E(u^+;G\times G)\le \mathcal E(u;G\times G)<\infty$ for any $G\Subset \Omega$.
Therefore $u^+\in H^s_{\rm loc}(\Omega)$, and $u^+\in X^s(\Omega;Z)$ follows. 

Next, assume that $\mathcal E(u;G\times G)= \mathcal E(u^+;G\times G)$ on some domain $G\Subset\Omega$. Then 
$$(u^-(x)-u^-(y))^2+2\Big(u^+(x)u^-(y)+u^-(x)u^+(y)\Big)=0$$ 
for a.e. $(x,y)\in G\times G$.
We infer that $u^-$ is constant a.e. on $G$. If $u^-=0$ then $u\ge 0$ on $G$;
{if $u^-\neq 0$ we get $u^+=0$, that is, $u\le 0$ on $G$.} The proof for the ''minus'' sign follows by replacing $u$ by $-u$.
\QED

\begin{Remark}
If $u\in L^1_{\rm loc}(U_1\cup U_2)$ and $\mathcal E_s(u;Z)$ is finite, then
$\mathcal E_s(|u|;Z)<\mathcal E_s(u;Z)$, {unless $u$ has constant sign on $U_1\cup U_2$.}
The proof runs with no changes.
\end{Remark}

{Our proof of Theorem} \ref{T:mp} requires the construction of a suitable barrier function.
The next Lemma slightly generalizes a result by {Ros-Oton and Serra  \cite{RoS}.}

\begin{Lemma}
\label{L:barrier}
Let $B_R(x^0)\subset \Omega$. For any $r\in(0, R)$ there exists a constant $c=c(R/r)>0$
and a function $\Phi\in H^s(\R^n)$
satisfying 
\begin{gather}
\nonumber
{\mathfrak L}^s_Z\Phi\le 0\qquad \text{in $B_R(x^0)\setminus \overline B_{r}(x^0)$}~;\\
\label{eq:Phi_data}
\Phi\equiv 1~~\text{in $B_{r}(x^0)$}~,\quad
\Phi\equiv 0~~\text{in $\R^n\setminus B_R(x^0)$}~,\quad
\Phi(x)\ge c(R-|x|)^s~~\text{in $B_R(x^0)$}~.
\end{gather}
\end{Lemma}

\proof
Without loss of generality we can assume $x^0=0$. 
Lemma 3.2 in \cite{RoS}, see also \cite[Lemma 2.2]{IMS}, provides the existence of $\Phi\in \widetilde H^s(B_R)$
satisfying (\ref{eq:Phi_data}) and $\Ds\Phi\le 0$ in $B_R\setminus \overline B_{{r}}$. To conclude
we claim that the distribution $\Ds\Phi-{\mathfrak L}^s_Z\Phi$ is nonnegative in $\Omega$.
Indeed, take $\eta\in {\cal C}^\infty_0(\Omega)$. Since both $\Phi$ and $\eta$ vanish on $\R^n\setminus\Omega$,
we have
\begin{eqnarray*}
\langle\Ds\Phi-{\mathfrak L}^s_Z\Phi,\eta\rangle&=&
\frac{C_{n,s}}{2}\iint\limits_{\R^{2n}\setminus Z}\frac{(\Phi(x)-\Phi(y))(\eta(x)-\eta(y))}{|x-y|^{n+2s}}~\!dxdy\\
&=&{C_{n,s}}\int\limits_{\Omega}\Phi(x)\eta(x)\Big(\int\limits_{\R^n\setminus (U_1\cup U_2)}\frac{dy}{|x-y|^{n+2s}}~\!dy\Big)dx,
\end{eqnarray*}
and the claim follows. In particular, ${\mathfrak L}^s_Z\Phi\le\Ds\Phi\le 0$ in $B_R\setminus \overline B_{{r}}$, and we are done.
\QED

\begin{Remark}
 It is worth to note that if $
 Z\subset Z'$ then for any nonnegative $\Phi\in\widetilde H^s(\Omega)$
the inequality ${\mathfrak L}^s_{Z}\Phi\le {\mathfrak L}^s_{Z'}\Phi$ holds in $\Omega$. The proof runs without changes.
\end{Remark}

We conclude this preliminary section by the following remark. We fix an exponent $\bp>2$; precisely we choose $\bp=4$ (for instance) if $n=1\le 2s$, and 
$\bp=\sstar=\frac{2n}{n-2s}$ if $n>2s$. Take any radius $r\in\big(1,2\big]$. The 
 Sobolev embedding theorem implies  
$\displaystyle{\mathcal E_s(u;\R^n\times\R^n) \ge 
c\,\Big(\int\limits_{B_r}|u|^\bp~\!dx\Big)^{\!2/\bp}}$ 
for any $u\in \widetilde H^s(B_r)$.

Now let $\rho\in\big(1,r\big)$. Since for $u\in \widetilde H^s(B_\rho)$ one has
$$
\mathcal E_s(u;\R^n\times\R^n)= \mathcal E_s(u;B_r\times B_r)+C_{n,s}\int\limits_{B_r}|u(x)|^2\,\Big(\int\limits_{\R^n\setminus B_r}\frac
{dy}{|x-y|^{n+2s}}\Big)dx,
$$
we plainly infer that 
\begin{equation}
\label{RSobolev}
\mathcal E_s(u;B_r\times B_r)+\frac{1}{(r-\rho)^{2s}}\int\limits_{B_r} |u|^2~\!dx\ge  c\,\Big(\int\limits_{B_\rho}|u|^\bp~\!dx\Big)^\frac2\bp
\quad \text{for any $u\in \widetilde H^s(B_\rho)$.}
\end{equation}

\section{Pointwise estimates for ${\mathfrak L}^s_Z$-subharmonic functions}

First, we prove a Caccioppoli-type inequality. We use again the notation introduced in (\ref{eq:3,1}).

\begin{Lemma}
\label{L:Cacc}
Let  $G\subseteq \Omega$ be a  {Lipschitz} domain, $w\in X^s(\Omega;Z)$ and $\f\in {\cal C}^\infty_0(G)$. Then
\begin{eqnarray}\label{eq:technicalA}
\mathcal E_s(\f w^+;G\times G)&\le&\langle {\mathfrak L}^s_Z w,\f^2 w^+\rangle\\
&+&\frac{C_{n,s}}{2}\iint\limits_{G\times G}w^+(x)w^+(y)\Psi_{\!\f}(x,y)~\!dxdy~\!
-\int\limits_G w^+\f^2(-\Delta_{(U_1\cup U_2)\setminus G}^{\!N}\!\!\left.\right)^s_{{\rm R}}\!w^+dx~\!.
\nonumber
\end{eqnarray}
\end{Lemma}

\proof
Note that all quantities in (\ref{eq:technicalA}) are finite by Lemmata \ref{L:dual}, \ref{L:3} and \ref{L:truncation}.
We compute 
\begin{eqnarray*}
&&(w(x)-w(y))((\f^2w^+)(x)-(\f^2w^+)(y))-((\f w^+)(x)-(\f w^+)(y))^2\\
&&\qquad\qquad=-w^+(x)w^+(y)(\f(x)-\f(y))^2+\big(\f(y)^2w^-(x)w^+(y)+\f(x)^2w^+(x)w^-(y)\big)\\
&&\qquad\qquad\ge -w^+(x)w^+(y)(\f(x)-\f(y))^2
\end{eqnarray*}
to infer 
\begin{eqnarray*}
\langle {\mathfrak L}^s_Z w,\f^2 w^+\rangle&=&
\frac{C_{n,s}}{2}\iint\limits_{Z}\frac{(w(x)-w(y))((\f^2w^+)(x)-(\f^2w^+)(y))^2}{|x-y|^{n+2s}}~\!\!dxdy\\
&\ge& \mathcal E_s(\f w^+;Z)-
\frac{C_{n,s}}{2}\iint\limits_{Z}w^+(x)w^+(y)\Psi_{\!\f} 
~\!\!dxdy\\
&=&\mathcal E_s(\f w^+,G\times G)+ \int\limits_G M^Z_G(x)|w^+\f|^2~\!dx
-\frac{C_{n,s}}{2}\iint\limits_{Z}w^+(x)w^+(y)\Psi_{\!\f} 
~\!\!dxdy
\end{eqnarray*}
by Lemma \ref{L:1}. {We compute} 
\begin{multline*}
\int\limits_G M^Z_G(x)|w^+\f|^2~\!dx=C_{n,s}\int\limits_Gw^+(x){|\f(x)|^2}
\Big(\int\limits_{{(U_1\cup U_2)\setminus G}}~\!\frac{(w^+(x)- w^+(y)+ w^+(y))}{|x-y|^{n+2s}}~\!dy\Big)dx
\\
\quad 
=\int\limits_Gw^+{|\f|^2}(-\Delta_{(U_1\cup U_2)\setminus G}^{\!N}\!\!\left.\right)^s_{{\rm R}}\!w^+dx+
C_{n,s}\int\limits_{G} w^+(x){|\f(x)|^2}~\Big(\!\!\!\! \int\limits_{(U_1\cup U_2)\setminus G}~\frac{w^+(y)}{|x-y|^{n+2s}}~\!dy\Big)dx~\!.
\end{multline*}
{Since $\Psi_\f\equiv 0$ on $\Gcomplement\times\Gcomplement$ we have, by (\ref{eq:sets}) with $\widetilde G=G$, }
\begin{multline*}
\frac{C_{n,s}}{2}\iint\limits_{Z}w^+(x)w^+(y)\Psi_{\!\f}~\!dxdy\\
=\frac{C_{n,s}}{2}\iint\limits_{G\times G}w^+(x)w^+(y)\Psi_{\!\f}~\!dxdy+
C_{n,s}\int\limits_{G}w^+(x)\Big(\!\!\!\!\!\int\limits_{(U_1\cup U_2)\setminus G}w^+(y)\Psi_{\!\f}~\!dy\Big)dx\\
=\frac{C_{n,s}}{2}\iint\limits_{G\times G}w^+(x)w^+(y)\Psi_{\!\f}~\!dxdy+
C_{n,s}\int\limits_{G}w^+(x)\f(x)^2\Big(\!\!\!\!\!\int\limits_{(U_1\cup U_2)\setminus G}\frac{w^+(y)}{|x-y|^{n+2s}}~\!dy\Big)dx~\!,
\end{multline*}
and the Lemma follows.
\QED

\begin{Remark}
Inequality (\ref{eq:technicalA}) was essentially proved in \cite[Theorem 1.4]{DKP},
in a weaker form but in a non-Hilbertian setting and for more general kernels. 
\end{Remark}

The next De Giorgi-type result is obtained by suitably modifying the argument for \cite[Theorem 1.1]{DKP}.

\begin{Lemma}
\label{L:DKP2} 
For any $u\in X^s(\Omega;Z)$ such that ${\mathfrak L}^s_Z u\le 0$~in $\Omega$ and for every ball $B_{2r}(x^0)\subseteq\Omega$, one has
\begin{equation}
\label{eq:DKP}
\sup_{B_r(x^0)} u
\le \Big(\frac {\widehat c}{r^n}\!\!\int\limits_{B_{2r}(x^0)}\!\!\!|u^+(x)|^2~\!dx\Big)^\frac12
+ 
r^{2s}\!\!\int\limits_{(U_1\cup U_2)\setminus B_{r}(x^0)}\frac{|u^+(x)|}{|x-x^0|^{n+2s}}~\!dx~\!,
\end{equation}
where $\widehat c>0$ depends only on $n$ and $s$. {In particular, $u$ is locally bounded from above
in $\Omega$.}
\end{Lemma}

%
%

\proof
First of all let us recall that {for $U_1\cup U_2=\R^n$} the last term in (\ref{eq:DKP}) 
is called {\em nonlocal tail}. For $Z\neq\R^n\times\R^n$ we call this term {\bf relative nonlocal tail} and denote it by $\text{Tail}_Z(u^+;x^0,r)$.

By rescaling we can assume without loss of generality that $r=1$ and $x^0=0$. 
{We introduce a parameter $\tk>0$ satisfying}
\begin{equation}
\label{eq:tk}
\tk\ge \text{Tail}_Z(u^+;0,1)
\end{equation}
(its value will be chosen later). For any integer $j\ge 0$
we put
$$
\begin{array}{lll}
r_j=1+2^{-j}~,& k_j=\tk(1-2^{-j}) ~,&
B_j=B_{r_j}~; \\
\tr_j=\displaystyle\frac{r_j+r_{j+1}}{2}~,&\tk_j=\displaystyle\frac{k_j+k_{j+1}}{2}{,} 
&\tB_j=B_{\tr_j}~;\\
w_j=(u-k_j)^+~,&\tw_j=(u-\tk_j)^+~,&
\alpha_j=\big(\displaystyle\int\limits_{B_j}|w_j|^2~\!dx\big)^\frac12.
\end{array}
$$
The following relations are obvious:
\begin{eqnarray}
\nonumber
&r_j\searrow1~,\quad r_{j+1}<\tr_j<r_{j}~;\quad 
&k_j\nearrow \tk~,\quad  k_j<\tk_j<k_{j+1};\\
&\tw_j\le w_j~,\quad\quad
&\tw_j\le \displaystyle\frac{w_j^2}{\tk_j-k_j}=\frac{ 2^{j+2}}{\tk}~\!w_j^2~;
\label{eq:tww}\\
&\alpha^2_0=\displaystyle\int\limits_{B_2}|u^+|^2~\!dx~,\quad\quad
&\alpha^2_j\to \int\limits_{B_1}|(u-\tk)^+|^2~\!dx\quad{\rm as}\quad j\to \infty.
\label{eq:alpha}
\end{eqnarray}
In addition, we have
\begin{equation}
\label{eq:tww0}
{w_{j+1}^2\Big(\frac{\tk}{2^{j+2}}\Big)^{\bp-2}=w_{j+1}^2(k_{j+1}-\tk_j)^{\bp-2}\le \tw_j^\bp}~\!,
\end{equation}
where the exponent $\bp>2$ was introduced at the {end of Section \ref{2}}.

Next, for any integer $j\ge 0$ we fix a cut-off function $\f_j$ satisfying
$$
\f_j\in {\cal C}^\infty_0(\tB_j)~,\quad 0\le\f_j\le 1~,\quad \f\equiv 1~~\text{on $B_{j+1}$},\quad
\|\nabla\f	\|_\infty\le 2^{j+3}.
$$
Since $\tw_j\f_j\in \widetilde H^s(B_j)$ and $1<\tr_j<r_j<2$, by (\ref{RSobolev}) {with
$\rho=\tr_j$ and $r=r_j$, we have}
\begin{equation}
c\Big(\int\limits_{B_j}|\tw_j\f_j|^\bp\Big)^\frac{2}{\bp}\le
\mathcal E_s(\tw_j\f_j;B_j\times B_j)+2^{2s(j+2)}\int\limits_{B_j}|\tw_j\f_j|^2~\!dx.
\label{eq:tildeSobolev}
\end{equation}

Notice that 
$\langle {\mathfrak L}^s_Z (u-\tk_j),\tw_j\f_j^2\rangle\le 0$, since
${\mathfrak L}^s_Z(u-\tk_j)= {\mathfrak L}^s_Z u\le0$ in $B_j$ 
and $\tw_j\f_j^2\in \widetilde H^s(B_j)$ {is nonnegative}. 
Using Lemma \ref{L:Cacc} with $w=u-\tk_j$, we 
infer
$$
\mathcal E_s(\tw_j\f_j;{B_j}\times {B_j})
\le 
c\iint\limits_{{B_j}\!\times\!{B_j}}\tw_j(x)\tw_j(y)\Psi_{\f_j}(x,y)~\!dxdy
-\int\limits_{B_j} \tw_j\f_j^2(-\Delta_{(U_1\cup U_2)\setminus B_j}^{\!N}\!\!\left.\right)^s_{{\rm R}}\!\tw_j\,dx,
$$
so that
\begin{equation}
\label{eq:rhs}
\Big(\int\limits_{B_j}|\tw_j\f_j|^\bp~\!dx\Big)^\frac{2}{\bp}\le c\Big(J_1-J_2+\tk^2~\!{2^{2sj}}~\!\Big(\frac{\alpha_{j}}{\tk}\Big)^2\Big)
\end{equation}
by (\ref{eq:tildeSobolev}), where
$$
J_1=\iint\limits_{{B_j}\!\times\!{B_j}}\tw_j(x)\tw_j(y)\Psi_{\f_j}(x,y)~\!dxdy~,\quad
J_2=\int\limits_{B_j} \tw_j\f_j^2(-\Delta_{(U_1\cup U_2)\setminus B_j}^{\!N}\!\!\left.\right)^s_{{\rm R}}\!\tw_j)\,dx.
$$
We estimate  from below the left-hand side of (\ref{eq:rhs}) {via} (\ref{eq:tww0}):
\begin{equation}
\label{eq:right}
\int\limits_{B_j}|\tw_j\f_j|^\bp~\!dx\ge 
\int\limits_{B_{j+1}}|\tw_j|^\bp~\!dx\ge c\,\Big(\frac{\tk}{2^{j}}\Big)^{\bp-2}\int\limits_{B_{j+1}}|w_{j+1}|^2~\!dx
=c~\!\tk^\bp~\!2^{j(2-\bp)}~\!\Big(\frac{\alpha_{j+1}}{\tk}\Big)^2.
\end{equation}
We estimate $J_1$ { by using}
$$
\Psi_{\f_j}(x,y)\le \|\nabla\f_j\|_\infty^2\,|x-y|^{-(n+2s-2)}
{\le c 2^{2j} \,|x-y|^{-(n+2s-2)}}
$$ 
and
the Cauchy-Bunyakovsky-Schwarz inequality, to obtain
\begin{eqnarray}
\nonumber
J_1&\le& c~\! 2^{2j}\iint\limits_{{B_j}\!\times\!{B_j}}\frac{\tw_j(x)}{|x-y|^{\frac{n+2s-2}{2}}}~\!
\frac{\tw_j(y)}{|x-y|^{\frac{n+2s-2}{2}}}~\!\!dxdy\le
 c~\! 2^{2j}\iint\limits_{{B_j}\!\times\!{B_j}}\frac{|\tw_j(x)|^2}{|x-y|^{n+2s-2}}~\!\!dxdy
\\
&=& c~\! 2^{2j}\int\limits_{B_j}|\tw_j(x)|^2\Big(\int\limits_{B_j}\frac{dy}{|x-y|^{n+2s-2}}~\!\Big)dx
{\,\le c r_j^{2-2s}~\!2^{2j} \alpha_j^2\le c~\!\tk^2~\!2^{2j}~\!\Big(\frac{\alpha_{j}}{\tk}\Big)^2.}
\label{eq:J1}
\end{eqnarray}
We handle $J_2$ as follows. For $x\in {\rm supp}(\f_j)\subset\tB_j$ 
and $y\in \Omega\setminus B_j$ we have 
$$
\frac{|y|}{|x-y|}\le 1+\frac{|x|}{|x-y|}\le 1+\frac{r_j}{r_j-\tr_j}\le c~\!2^{j}.
$$
Hence, using also (\ref{eq:tww}) we can estimate
$$
\tw_j(x)|\f_j(x)|^2\,\frac{\tw_j(y)-\tw_j(x)}{|x-y|^{n+2s}}\le \frac{c}{\tk}~\!|w_j(x)|^2~\!{2^{j(n+2s+1)}}~\!\frac{w_j(y)}{|y|^{n+2s}}~\!,
$$
so that
\begin{eqnarray*}
-J_2&=&\int\limits_{B_j} \tw_j(x)|\f_j(x)|^2\Big(\int\limits_{{(U_1\cup U_2)\setminus {B_j}}}\frac{\tw_j(y)-\tw_j(x)}{|x-y|^{n+2s}}dy\Big)~\!dx\\
&\le&\frac{c}{\tk}~\!2^{j(n+2s+1)}~\!\Big(\int\limits_{{(U_1\cup U_2)\setminus {B_j}}}\frac{w_j(y)}{|y|^{n+2s}}dy\Big)~\!
\int\limits_{B_j} {|w_j|^2}~\!dx\\
&\le&c \tk~\!{2^{j(n+2s+1)}}~\text{Tail}_Z(u^+;0,1)~\!\Big(\frac{\alpha_{j}}{\tk}\Big)^2
\end{eqnarray*}
because $B_j\supset B_1$ and $w_j\le u^+$.
Comparing with (\ref{eq:rhs}), (\ref{eq:right}) and (\ref{eq:J1}) we arrive at
$$
2^{\frac{2(2-\bp)}{\bp}j}~\!\Big(\frac{\alpha_{j+1}}{\tk}\Big)^\frac4\bp \le
c~\!{2^{j(n+2s+1)}}\Big(1+{\tk}^{-1}~\!\text{Tail}_Z(u^+;0,1)\Big)\Big(\frac{\alpha_{j}}{\tk}\Big)^2~\!.
$$
Taking (\ref{eq:tk}) into account, we can conclude that
\begin{equation}
\label{eq:ind}
\frac{\alpha_{j+1}}{\tk}\le (\widehat c^{\frac \beta2}\eta^{-\frac1\beta})\eta^j
\Big(\frac{\alpha_{j}}{\tk}\Big)^{\beta+1},
\end{equation}
where
{$\displaystyle{\beta=\frac\bp2-1>0~,~~ \eta=2^{\frac{\bp}{4}(n+2s+1)+\beta}>1}$.}
Now we choose the free parameter $\tk$, namely
$$
\tk= \text{Tail}_Z(u^+;0,1)+\widehat c^{\frac 12}\alpha_0=
\text{Tail}_Z(u^+;0,1)+\Big(\widehat c\int\limits_{B_1}|u^+|^2~\!dx\Big)^\frac12,
$$
compare with (\ref{eq:DKP}). The above choice of $\tk$ guarantees that 
\begin{equation}
\label{eq:induction}
\widehat c^{\frac 12}~\!\frac{\alpha_j}{\tk}\le\eta^{-\frac{j}{\beta}}
\end{equation}
for $j=0$. Using induction and (\ref{eq:ind}) one easily gets that (\ref{eq:induction}) holds for any $j\ge 0$. Thus
 $\alpha_j\to 0$ and hence $(u-\tk)^+\equiv 0$ on $B_1$ by (\ref{eq:alpha}). The proof is complete.
\QED

\section{Main results} 
\label{S:main}

We are in position {to state and prove} a strong maximum principle for the nonlocal operator ${\mathfrak L}^s_Z$,  that is the main result of the present paper. 
\begin{Theorem}
\label{T:mp}
Let $u$ be a nonconstant measurable function on $U_1\cup U_2$ such that 
$$
u\in {H^s_{\rm loc}}(\Omega)~,\quad \int\limits_{U_1\cup U_2}\frac{|u(x)|}{1+|x|^{n+2s}}~\!dx<\infty~,\quad {\mathfrak L}^s_Z u\ge 0~~\text{in $\Omega$.}
$$
Then $u$ is lower semicontinuous on $\Omega$, locally bounded from below on $\Omega$ and
$$
\displaystyle{u(x)>\inf_{U_1\cup U_2} u}\qquad\text{for every $x\in \Omega$.}
$$ 
\end{Theorem}

\proof
{First, local boundedness from below follows from Lemma \ref{L:DKP2}.}

To check the first claim it suffices to show that $u$ has a representative that is lower semicontinuous
on any fixed domain $G\Subset \Omega$. From 
$$
\frac{C_{n,s}}{2}\int\limits_G\Big(\int\limits_{G} \frac{(u(x)-u(y))^2}{|x-y|^{n+2s}}~\!dy\Big)dx<\infty
$$ 
we infer that
\begin{equation}
\label{eq:a.e.}
\int\limits_{G} \frac{(u(x)-u(y))^2}{|y-x|^{n+2s}}~dy<\infty
\end{equation}
for a.e. $x\in G$. Let $G_0$ be the set of Lebesgue points $x\in G$ for $u$ that satisfy (\ref{eq:a.e.}).
We can assume that
$u(x^0)=\displaystyle\liminf_{x\to x^0} u(x)$ for any $x^0\in G\setminus G_0$, because $ G\setminus G_0$
has null Lebesgue measure.

Our next goal is to show that $u(x^0)\le\displaystyle\liminf_{x\to x^0} u(x)$ for any $x^0\in G_0$. We use 
Lemma \ref{L:DKP2} with $u$ replaced by $u(x^0)-u$ to get
\begin{equation}
\label{eq:forni}
\inf_{B_{r}(x^0)} u
\ge u(x^0)-\text{\rm Tail}_{Z}((u(x^0)-u)^+;x^0,r)-\Big(\frac{\widehat c}{r^n}\!\!\int\limits_{~B_{2r}(x^0)}\!\!\!\!|(u(x^0)-u)^+|^2~\!dx\Big)^\frac12
\end{equation}
for any $r>0$ small enough. First we split
\begin{multline*}
\text{\rm Tail}_{Z}((u(x^0)-u)^+;x^0,r)=
r^{2s}\!\int\limits_{(U_1\cup U_2)\setminus B_{2r}(x^0)}\frac{|u(x^0)-u(x)|}{|x-x^0|^{n+2s}}~\!dx\\
\le 
r^{2s}\!\int\limits_{(U_1\cup U_2)\setminus G}\frac{|u(x^0)|+|u(x)|}{|x-x^0|^{n+2s}}~\!dx
~+ ~
r^{2s}\!\int\limits_{G\setminus B_{2r}(x^0)}\frac{|u(x^0)-u(x)|}{|x-x^0|^{n+2s}}~\!dx
=:P_r+Q_r.
\end{multline*}
{We readily obtain}
$$
{P_r} \le c(\text{dist}(x^0,\partial G))~\!r^{2s}\int\limits_{(U_1\cup U_2)}\frac{|u(x^0)|+|u(x)|}{1+|x|^{n+2s}}~\!dx
\to 0
$$
as $r\to 0$. Next we use the Cauchy-Bunyakovsky-Schwarz inequality
to estimate 
$$
{Q_r\le  r^{2s}\Big(\int\limits_{G}\frac{(u(x^0)-u(x))^2}{|x-x^0|^{n+2s}}~\!dx\Big)^{\!\frac12}\!\Big(\int\limits_{\R^n\setminus B_{2r}(x^0)}\frac{dx}{|x-x^0|^{n+2s}}\Big)^{\frac12}\!=
c r^s \Big(\int\limits_{G}\frac{(u(x^0)-u(x))^2}{|x-x^0|^{n+2s}}~\!dx\Big)^{\frac12}\!.}
$$
Since (\ref{eq:a.e.}) is satisfied at $x=x^0$, we have that  
{$Q_r\to 0$. Thus $\text{\rm Tail}_{Z}((u(x^0)-u)^+;x^0,r)\to 0$ as $r\to 0$. }
Further, the last term in (\ref{eq:forni}) goes to zero as $r\to 0$ because $x^0$ is a Lebesgue point for $u$. 
Thus $\displaystyle\liminf_{x\to x^0} u(x)\ge u(x^0)$, and the first statement is proved.

{Next, assume by contradiction} that $u$ is bounded from below and
$$
\Omega_+:=\{x\in \Omega~|~ u(x)>m:=\inf\limits_{U_1\cup U_2} u~\}
$$
is strictly contained in $\Omega$. Since $u$ is lower semicontinuous on ${\Omega}$, {the set ${\Omega_+}$ is open} and has a nonempty
boundary in ${\Omega}$. 

Fix a point $\xi\in{\Omega}\cap \partial {\Omega_+}$, so that $u(\xi)=m$. Using again the lower-semicontinuity
of $u$, we can find $R>r>0$ and a point $x^0\in {\Omega_+}$, such that
$\xi\in B_R(x^0)\Subset {\Omega}$ and $u(x)\ge \frac12(u(x^0)+m)>m$ for every $x\in B_{{r}}(x^0)$.
We can assume that $x^0=0$ to simplify notations.
Thus we have the following situation: 
\begin{equation}
\label{eq:xi}
\xi\in B_R\subset{\Omega}~,\quad u(\xi)=m
~,\quad \inf_{B_{{r}}}u(x) \ge m+\delta
\end{equation}
for some $\delta>0$.
Let $\Phi$ be the function defined in Lemma \ref{L:barrier}. We claim that $u\ge m+\delta\Phi>m$ in 
$B_R\setminus \overline B_{{r}}$, that gives a contradiction
with (\ref{eq:xi}). 

Indeed, define $v=u-\delta\Phi$, so that 
$$
v= u-\delta\ge m~~\text{in $B_{{r}}$},\quad 
v=u\ge m~~\text{in $(U_1\cup U_2)\setminus B_R$.}
$$
Our goal is to show that $v\ge m$ also on $B_R\setminus \overline B_{{r}}$.

Clearly $v\in X^s({\Omega};Z)$ as $u, \Phi \in X^s({\Omega};Z)$. By Lemma \ref{L:truncation} 
this implies $v^m_\pm:=(v-m)^\pm\in X^s({\Omega};Z)$. Next, notice that $v^m_-=0$ out of $\overline B_R\setminus {B}_{{r}}$. 
Therefore $v^m_-\in\widetilde H^s(B_R\setminus \overline B_{{r}})$, and using Lemma \ref{L:dual} we obtain
\begin{equation}
\label{eq:god}
\langle {\mathfrak L}^s_Z v, v^m_-\rangle = \langle {\mathfrak L}^s_Z u, v^m_-\rangle-\delta \langle {\mathfrak L}^s_Z \Phi, v^m_-\rangle
\ge 0.
\end{equation}
However,
\begin{eqnarray*}
\langle {\mathfrak L}^s_Z v, v^m_-\rangle
&=&\frac{C_{n,s}}{2}\iint\limits_{\!\! Z}\frac{((v(x)-m)-(v(y)-m))(v^m_-(x)-v^m_-(y))}{|x-y|^{n+2s}}~\!dxdy\\
&=& -\frac{C_{n,s}}{2}\iint\limits_{\!\! Z}\frac{v^m_+(x)v^m_-(y)+v^m_+(y)v^m_-(x)}{|x-y|^{n+2s}}~\!dxdy-\mathcal E_s(v^m_-;Z)\\
&\le& -\mathcal E_s(v^m_-;B_R\times B_R),
\end{eqnarray*}
so that (\ref{eq:god}) implies $\mathcal E_s(v^m_-;B_R\times B_R)=0$, that together with 
$v^m_-\in\widetilde H^s(B_R\setminus \overline B_{{r}})$ gives the conclusion.
\QED

By choosing $U_1=U_2=\R^n$ we have $Z=\R^n\times\R^n$. It is well known that
${\mathfrak L}^s_{\R^n\times\R^n}u=\Ds u$ pointwise on $\R^n$ if $u\in C^2(\R^n)$. Thanks to Lemma \ref{L:dual}, we can say that ${\mathfrak L}^s_{\R^n\times\R^n}u=\Ds u$ in
$\widetilde H^s(G)'$, for every $u\in X^s(\Omega;\R^n\times\R^n)$ and for any Lipschitz domain $G\Subset \Omega$. 
From Theorem \ref{T:mp} we immediately infer the next result. 

\begin{Corollary}[Dirichlet Laplacian]
\label{C:Dir}
Let $\Omega\subseteq\R^n$ be a domain, and let $u$ be a nonconstant measurable function on $\R^n$ such that 
$$
u\in {H^s_{\rm loc}}(\Omega)~,\quad \int\limits_{\R^n}\frac{|u(x)|}{1+|x|^{n+2s}}~\!dx<\infty~,\quad \Ds u\ge 0~~\text{in $\Omega$.}
$$
Then $u$ is lower semicontinuous on $\Omega$, locally bounded from below on $\Omega$ and 
 $\displaystyle{u(x)>\inf_{\R^n} u}$ for every $x\in \Omega$.
\end{Corollary}

The Restricted Laplacian is obtained by choosing $U_1=U_2=\Omega$. In fact, 
${\mathfrak L}^s_{\Omega\times\Omega}u=\DsNROmega u$  on $\Omega$ if $u\in C^2(\Omega)$,
and ${\mathfrak L}^s_{\Omega\times\Omega}u=\DsNROmega u$ in
$\widetilde H^s(G)'$, for every $u\in X^s(\Omega;\R^n\times\R^n)$ and every Lipschitz domain $G\Subset \Omega$. 
From Theorem \ref{T:mp} we infer the next result.

\begin{Corollary}[Restricted Laplacian]
\label{C:Res}
Let $\Omega\subset\R^n$ be a domain, and let $u$ be a nonconstant measurable function on $\Omega$ such that 
$$
u\in {H^s_{\rm loc}}(\Omega)~,\quad \int\limits_{\Omega}\frac{|u(x)|}{1+|x|^{n+2s}}~\!dx<\infty~,\quad 
\DsNROmega u\ge 0~~\text{in $\Omega$.}
$$
Then $u$ is lower semicontinuous on $\Omega$, locally bounded from below on $\Omega$ and
$\displaystyle{u(x)>\inf_{\Omega} u}$~~for every $x\in \Omega$. \end{Corollary}

Next, we choose $U_1=\Omega, U_2=\R^n$, so that $Z=\R^{2n}\setminus(\complement)^2$. 
By \cite[Lemma 3]{DRoV} we have that ${\mathfrak L}^s_{\R^{2n}\setminus(\complement)^2} u=\DsNSOmega u$
 if $u\in C^2(\R^n)\cap L^\infty(\R^n)$,  compare with (\ref{eq:SR}).
From the computations there and {thanks to} Lemma \ref{L:dual} we can identify the distributions
${\mathfrak L}^s_{\R^{2n}\setminus(\complement)^2}u$ and $\DsNSOmega u$  for functions
$u\in X^s(\Omega;\R^{2n}\setminus(\complement)^2)$ (see also  \cite[Definition 3.6]{DRoV}). 
Theorem \ref{T:mp} immediately implies
the next corollary, see also \cite[Theorem 1.1]{BM} for a related result.

\begin{Corollary}[Semirestricted Laplacian]
\label{C:Sem}
Let $\Omega\subset\R^n$ be a domain, and let $u$ be a nonconstant measurable function on $\R^n$ such that 
$$
u\in {H^s_{\rm loc}}(\Omega)~,\quad \int\limits_{\R^n}\frac{|u(x)|}{1+|x|^{n+2s}}~\!dx<\infty~,\quad \DsNSOmega u\ge 0~~\text{in $\Omega$.}
$$
Then $u$ is lower semicontinuous on $\Omega$, locally bounded from below on $\Omega$ and
 $\displaystyle{u(x)>\inf_{\R^n} u}$ for every $x\in \Omega$.
\end{Corollary}

We conclude by recalling that in the local case $s=1$, 
the strong maximum principle states that every nonconstant superharmonic function $u$ on $\Omega$
satisfies $\displaystyle{u(x)>\inf_{\Omega} u}$ for every $x\in \Omega$. 
Notice that in the non local, Neumann Restricted case we reached the same conclusion. In contrast, 
in the Dirichlet and in the Semirestricted cases a similar result can not hold, see the example in the next remark.

\begin{Remark}
\label{R:Alex}
Take any bounded domain $\Omega\in \R^n$ and two nonnegative functions $u,\psi\in {\cal C}^\infty_0(\R^n)$ such that
$0\le u\le 1$, $u\equiv 1$ on $\Omega$, $\text{supp}~\!\psi\subset\Omega$. For any $x\in \Omega$ we have
$$
\Ds u(x)= \DsNSOmega u(x)={C_{n,s}}~\!\cdot~\!{\rm P.\!V.}\irn\frac{1-u(y)}{|x-y|^{n+2s}}~dy>0.
$$
Since $\Ds u, \Ds\psi$ are smooth functions, we have that $\Ds(u-\eps\psi)\ge 0$ in $\Omega$, for some small
$\eps>0$. Then $u-\eps\psi$ satisfies the assumptions in Corollaries \ref{C:Dir} and \ref{C:Sem}, but 
$\displaystyle{\inf_{\Omega} (u-\eps\psi)}$ is achieved in $\Omega$. Clearly, $\displaystyle{\inf_{\R^n} (u-\eps\psi)}=0$
is not achieved in $\Omega$.
\end{Remark}

\appendix

\section{\!\!\!\!\!\!ppendix}


{We start with a proposition in fact proved in \cite{Sil}. It gives
the same conclusion as in Corollary \ref{C:Dir} under 
weaker summability assumptions on $u$. Notice however that
$n>2s$ is needed (this is a restriction only if $n=1$), and that Silvestre's 
construction cannot be easily extended to more general operators
such as the Restricted and Semirestricted ones.
\begin{Proposition}
\label{P:Sil}
Assume $n>2s$ and let $u$ be a nonconstant measurable function on $\R^n$
such that $\frac{u(x)}{1+|x|^{n+2s}}\in L^1(\R^n)$ and $\Ds u\ge 0$ in the distributional
sense on $\Omega$, that is,
$$
\langle \Ds u,\f\rangle=\irn u\Ds \f~\!dx \ge 0 \qquad\text{for any~
$\f\in{\cal C}^\infty_0(\R^n), ~~\f\ge 0$.}
$$
Then $u$ is lower semicontinuous on $\Omega$ and  $\displaystyle{u(x)>\inf_{\R^n} u}$ for every $x\in \Omega$.
\end{Proposition}
\proof
First, notice that $\Ds u$ is a well defined distribution, as
$(1+|x|^{n+2s})\Ds\f$ is a bounded function on $\R^n$, for 
any $\f\in C^\infty_0(\R^n)$. 
{Proposition 2.2.6 in \cite{Sil} gives the lower semicontinuity of
$u$ in $\Omega$ and the relations
$u(x_0)\ge \displaystyle{\int_{\R^n} u(x)\gamma_r^s(x_0-x)~\!dx}>-\infty$
for any ball $B_r(x_0)\subset\Omega$,
where $\gamma_r^s$ is certain continuous and positive function on $\R^n$.
If $u$ is unbounded from below we are done; otherwise, 
we can assume $\displaystyle\inf_{\R^n} u=0$. Suppose that there exists $x_0\in\Omega$ such that 
$u(x_0)=0$. Take a ball $B_r(x_0)\subset\Omega$. Then
$0\ge \displaystyle{\int_{\R^n} u(x)\gamma_r^s(x_0-x)~\!dx}\ge 0$, that immediately implies $u\equiv 0$ in $\R^n$, a contradiction.}
%
\QED
}

Now we recall that the Spectral Dirichlet/Neumann fractional Laplacian is the $s$-th power of standard Dirichlet/Neumann Laplacian 
in $\Omega$ in the sense of spectral theory.

A strong maximum principle for the Spectral Dirichlet Laplacian follows from \cite[Lemma 2.6]{CDDS} 
and reads as follows.

\begin{Proposition}
 Let \ $\Omega\subset\R^n$ be a bounded domain, and let a function $u\in \Hs$ be such that $\DsSpOmega u\ge0$ in $\Omega$ the sense of distributions. Then either $u\equiv0$ or 
 $\inf\limits_Ku>0$ for arbitrary compact set $K\subset \Omega$.
\end{Proposition}

A strong maximum principle for the Spectral Neumann Laplacian can be  obtained from the results 
in \cite{CDDS}.

\begin{Theorem}
Let \ $\Omega\subset\R^n$ be a bounded Lipschitz domain, and let a function $u\in H^s(\Omega)$ be such that $\DsNSpOmega u\ge0$ in a subdomain $G\subset\Omega$ in the sense of 
distributions. Then either $u\equiv {\rm const}$ or $\inf\limits_Ku>\inf\limits_\Omega u$  for arbitrary compact set $K\subset G$.
 \end{Theorem}

\proof
It is well known, see \cite{ST}, \cite{CaSt} and \cite{AEW} for a general setting, that for any $u\in H^s(\Omega)$ the boundary value problem
\begin{equation}
-\div (y^{1-2s}\nabla w)=0\quad \mbox{in}\quad \Omega\times\mathbb R_+;\qquad w\big|_{y=0}=u;
\qquad \partial_{\bf n}w\big|_{x\in\partial\Omega}=0,
\label{eq:ST}
\end{equation}
has a unique weak solution $w_s^N(x,y)$, and 
\begin{equation*}
\DsNSpOmega u(x)=-\,\frac{2^{2s-1} \Gamma(s)}{\Gamma(1-s)}\cdot\lim\limits_{y\to0^+} 
y^{1-2s}\partial_yw_s^N(x,y)
\end{equation*}
(the limit is understood in the sense of distributions).

Without loss of generality we can assume that $\inf\limits_\Omega u=0$. Then by the maximum principle for (\ref{eq:ST}) we have $w\ge0$. By \cite[Lemma 2.6]{CDDS} the 
statement follows.
\QED

\bigskip
\medskip
\noindent
{\bf Acknowledgments}.
This paper was completed while the second author was visiting SISSA (Trieste), in 
January 2017. The authors wish to thank SISSA for the hospitality.

\footnotesize
\label{References}


\begin{thebibliography}{XX}
\footnotesize


\bibitem{AEW}
W. Arendt, A. F. M. ter Elst\ and\ M. Warma, {Fractional powers of sectorial operators via the Dirichlet-to-Neumann operator}, preprint arXiv:1608.05707 (2016).

\bibitem{BM}
B. Barrios and M. Medina, Strong maximum principles for fractional elliptic and parabolic
problems with mixed boundary conditions, preprint arXiv:1607.01505v2  (2016).

\bibitem{CaSt}
L. A. Caffarelli\ and\ P. R. Stinga, Fractional elliptic equations, Caccioppoli estimates and regularity, 
Ann. Inst. H. Poincar\'e Anal. Non Lin\'eaire {\bf 33} (2016), no.~3, 767--807.

\bibitem{CDDS}
A. Capella, J. D\'avila, L. Dupaigne\ and\ Y. Sire, Regularity of radial extremal solutions for some non-local semilinear equations, Comm. Partial Differential Equations {\bf 36} (2011), no.~8, 1353--1384.

\bibitem{DQ}
L. M. Del Pezzo\ and\ A. Quaas, A Hopf's lemma and a strong minimum principle for
the fractional $p$-Laplacian, preprint arXiv:1609.04725 (2016).

\bibitem{DKP}
A. Di Castro, T. Kuusi\ and\ G. Palatucci, Local behavior of fractional $p$-minimizers, 
Ann. Inst. H. Poincar\'e Anal. Non Lin\'eaire {\bf 33} (2016), no.~5, 1279--1299.

\bibitem{DRoV}
S. Dipierro, X. Ros-Oton\ and\ E. Valdinoci,
Nonlocal problems with Neumann boundary conditions, Rev. Mat. Iberoam., to appear.

\bibitem{G}
G. Grubb, Regularity of spectral fractional Dirichlet and Neumann problems, Math. Nachr. {\bf 289} (2016), no.~7, 831--844. 

\bibitem{G1}
Q. Y. Guan, Integration by parts formula for regional fractional Laplacian, Comm. Math. Phys. {\bf 266} (2006), no.~2, 289--329. 

\bibitem{G2}
Q. Y. Guan and Z. M. Ma, Boundary problems for fractional Laplacians, Stoch. Dyn., 5
(2005), 385--424.

\bibitem{IMS}
A. Iannizzotto, S. Mosconi\ and\ M. Squassina, 
$H\sp s$ versus $C\sp 0$-weighted minimizers, NoDEA Nonlinear Differential Equations Appl. {\bf 22} (2015), no.~3, 477--497.

\bibitem{FL} R. Musina\ and\ A. I. Nazarov, On fractional Laplacians, Comm. Partial Differential Equations 
{\bf 39} (2014), no.~9, 1780--1790. 

\bibitem{Ro}
X. Ros-Oton, Nonlocal elliptic equations in bounded domains: a survey, Publ. Mat. {\bf 60} (2016), no.~1, 3--26.

\bibitem{RoS}
X. Ros-Oton\ and\ J. Serra, The Dirichlet problem for the fractional Laplacian: regularity up to the boundary, 
J. Math. Pures Appl. (9) {\bf 101} (2014), no.~3, 275--302.

\bibitem{Sil}
L. Silvestre, Regularity of the obstacle problem for a fractional power of the Laplace operator, Comm. Pure Appl. Math. {\bf 60} (2007), no.~1, 67--112.

\bibitem{ST}
P. R. Stinga\ and\ J. L. Torrea, Extension problem and Harnack's inequality for some fractional operators, 
Comm. Partial Differential Equations {\bf 35} (2010), no.~11, 2092--2122.

\bibitem{Tr}
H. Triebel, {\it Interpolation theory, function spaces, differential operators}, 
Deutscher Verlag Wissensch., Berlin, 1978. 

\bibitem{W1}
M. Warma, The $p$-Laplace operator with the nonlocal Robin boundary conditions on arbitrary open sets, 
Ann. Mat. Pura Appl. (4) {\bf 193} (2014), no.~1, 203--235.

\bibitem{W2}
M. Warma, The fractional relative capacity and the fractional Laplacian 
with Neumann and Robin boundary conditions on open sets, Potential Anal. {\bf 42} (2015), no.~2, 499--547. 

\bibitem{W3}
M. Warma, A fractional Dirichlet-to-Neumann operator on bounded Lipschitz domains, 
Commun. Pure Appl. Anal. {\bf 14} (2015), no.~5, 2043--2067.

\end{thebibliography}
\end{document}